\documentclass{amsart}
\usepackage{amssymb,amsthm}
\usepackage{a4}
\makeatletter
  
  \@addtoreset{equation}{section}
 \makeatother

\newtheorem{theorem}{Theorem}[section]
\newtheorem{lemma}[theorem]{Lemma}
\newtheorem{corollary}[theorem]{Corollary}
\newtheorem{remark}[theorem]{Remark}

\newtheorem{de}[theorem]{Definition}

\newcommand{\g}[2]{\ensuremath{\langle #1,#2 \rangle}}%metricf
\newcommand{\ja}{\ensuremath{\mathcal{J}_R}}%general Jacobi operator
\newcommand{\jp}{\ensuremath{\mathcal{J}_R(\pi_x)}}%Jacobi operator in \pi_x direction
\newcommand{\tjp}{\ensuremath{\mathcal{\tilde J}_R(\pi_x)}}%secondary Jacobi operator in \pi_x direction
\newcommand{\f}{\ensuremath{\mathcal{F}}}%Clifford family F

\begin{document}

\title{COMPLEX OSSERMAN ALGEBRAIC CURVATURE TENSORS AND CLIFFORD
FAMILIES}
\author{M. Brozos-V{\'a}zquez and P. Gilkey}
\begin{address}{MBV:Department of Geometry and Topology, Faculty of Mathematics, University of
Santiago de Compostela, 15782 Santiago de Compostela,
Spain}\end{address}
\begin{email}{mbrozos@usc.es}\end{email}
\begin{address}{PG: Mathematics Department, University of Oregon, Eugene, OR 97403, USA}\end{address}
\begin{email}{gilkey@uoregon.edu}\end{email}
\begin{abstract}
We use methods of algebraic topology to study the eigenvalue structure
of a complex Osserman algebraic curvature tensor. We classify the algebraic
curvature tensors which are both Osserman and complex Osserman in all but
a finite number of exceptional dimensions.
\end{abstract}
\maketitle

\section{Introduction}

In recent years, the Osserman problem has played an important role
in the understanding of curvature. The real setting has been studied previously; in this paper, we study
the complex setting. We introduce the following notational conventions. Let $\nabla$ be the Levi-Civita
connection of a Riemannian manifold $(M,g)$ and let $R$ be the associated Riemann curvature tensor:
$$R(x,y,z,w):=g((\nabla_x\nabla_y-\nabla_y\nabla_x-\nabla_{[x,y]})z,w)\,.$$
The
Jacobi operator
$\mathcal{J}_R$ and the skew-symmetric curvature operator $\mathcal{R}$ are characterized
by the identities:
\begin{equation}\label{eqn-1.a}
  g(\mathcal{J}_R(x)y,z)=R(y,x,x,z)\quad\text{and}\quad
  g(\mathcal{R}(x,y)z,w)=R(x,y,z,w)\,.
\end{equation}
Motivated by the seminal paper of Osserman
\cite{Osserman}, one says that $(M,g)$ is {\it Osserman} if the eigenvalues of
$\mathcal{J}_R$ are constant on the sphere bundle $S(M,g)$ of unit tangent vectors.
Since the local isometries of a local two-point homogeneous manifold act transitively on
$S(M,g)$, such manifolds are Osserman. Osserman wondered if the converse was also
true, that is, are Osserman manifolds necessarily local two-point-homogeneous
spaces. This question has been called the Osserman conjecture by subsequent authors
and has been also considered in the pseudo-Riemannian context; in this paper, we
will only work in the Riemannian context and refer to
\cite{E-K-R,Gilkey-book} for a discussion of the pseudo-Riemannian setting.

\subsection{Algebraic curvature tensors}\label{sect-1.1}
It turned out to be convenient
to work in a purely algebraic context in studying the Osserman conjecture. Let
$\mathfrak{M}:=(V,\langle\cdot,\cdot\rangle,R)$ be a {\it model}. This means that
$V$ is a vector space of dimension
$n$ which is equipped with a positive definite inner product
$\g\cdot\cdot$ and that
$R\in\otimes^4V^*$ is an {\it algebraic curvature tensor}, i.e. $R$ satisfies the
Riemannian curvature tensor identities:
\[
\begin{array}{l}
R(x,y,z,t)=-R(y,x,z,t)=R(z,t,x,y),\\
R(x,y,z,t)+R(y,z,x,t)+R(z,x,y,t)=0.
\end{array}
\]
One uses Equation (\ref{eqn-1.a}) to define the Jacobi operator
and skew-symmetric curvature operator in this setting as well;
$\mathfrak{M}$ is said to be {\it Osserman} if the eigenvalues of
$\ja(\cdot)$ are constant on the sphere
$S(V,\langle\cdot,\cdot\rangle)$ of unit vectors in $V$. Clearly,
if $(M,g)$ is a Riemannian manifold, and if $P\in M$, then
$\mathfrak{M}_P:=(T_PM,g_P,R_P)$ defines a model. Conversely,
every model is geometrically realizable.

\subsection{The Osserman conjecture}\label{sect-1.2}
This conjecture for Riemannian manifolds was established by Chi
\cite{Chi} in dimensions $n\equiv1$ (mod $2$), $n\equiv 2$ (mod
$4$) and $n=4$. Subsequent work by Nikolayevsky
\cite{nikolayevsky-1,nikolayevsky-2} has established the Osserman
conjecture in dimensions $n\ne16$; the case $n=16$ is still open.
Nikolayevsky used a two step approach following the discussion in
\cite{Gilkey-Swann-Vanhecke}. He first showed that any Osserman
model is given by a Clifford family as specified in
Equation (\ref{eqn-1.e}) of Section \ref{sect-1.7} below except in
dimension $16$. He then used the Second Bianchi Identity to prove
the necessary integrability results to show any Osserman manifold
of dimension $n\ne16$ was locally isometric to a rank $1$ symmetric space or was flat. Note that the
algebraic classification fails if $n=16$; indeed the curvature tensor of the Cayley plane is Osserman but
it is not given by a Clifford family, i.e. it is not expressible in the
form given in Equation (\ref{eqn-1.e}).

\subsection{The higher order Jacobi operator}\label{sect-1.3}
There are other related
questions. One may follow the discussion of Stanilov and Videv
\cite{StanilovVidev} to define a {\it higher order Jacobi operator} as follows. Let
$\{e_1,...,e_p\}$ be an orthonormal basis for a
$p$-plane $\mathcal{P}$.  Set
\[
\ja(\mathcal{P})=\sum_{i=1}^p \ja (e_i)\,;
\]
this is independent of the particular orthonormal basis chosen. If $p=1$, one recovers the
ordinary Jacobi operator. 
Furthermore, if $p=n$, then $\rho:=\mathcal{J}_R(V)$ is the {\it Ricci operator}; thus
the higher order Jacobi operator can also be thought of as a generalization of the
Ricci operator to lower dimensional subspaces. One says that a model
$\mathfrak{M}$ is {\it
$p$-Osserman} if the eigenvalues of
$\ja(\mathcal{P})$ are constant on the Grassmannian $Gr_p(V)$ of $p$-planes. The
geometry is very rigid in this setting. If $p=1$ or if $p=n-1$, then 
$\mathfrak{M}$ is $p$-Osserman if and only if $\mathfrak{M}$ is Osserman.
Thus these values of $p$ may be excluded from consideration.
If $2\le p\le n-2$, then it is known
\cite{gilkey-p-Osserman} that $\mathfrak{M}$ is
$p$-Osserman if and only if
$\mathfrak{M}$ has constant sectional curvature $c$, i.e. that $R=cR_0$ where
$R_0$ is given by:
\begin{equation}\label{eqn-1.b}
R_0(x,y,z,t):=\langle x,t\rangle\langle y,z\rangle-\langle
x,z\rangle\langle y,t\rangle\,.
\end{equation}

\subsection{Complex geometry}\label{sect-1.4}
In this paper, we will consider a complex analogue of these
questions.  Let $J$ denote an Hermitian almost complex structure
on $(V,\g\cdot\cdot)$; this means that $J$ is an isometry of
$(V,\langle\cdot,\cdot\rangle)$ with $J^2=-\operatorname{id}$. A
$2$-plane is said to be holomorphic if it is $J$-invariant and a
real linear transformation $T$ of $V$ is said to be complex linear
if $TJ=JT$. We let $\mathbb{CP}(V,\langle\cdot,\cdot\rangle,J)$ be
the set of all holomorphic $2$ planes. If $x\in
S(V,\langle\cdot,\cdot\rangle)$ is a unit vector, let
$\pi_x:=\operatorname{Span}\{x,Jx\}$. The natural map
$x\rightarrow\pi_x$ defines the Hopf fibration from
$S(V,\langle\cdot,\cdot\rangle)$ to
$\mathbb{CP}(V,\langle\cdot,\cdot\rangle,J)$. Let
$$
\ja(\pi_x):=\ja(x)+\ja(Jx)
$$
be the complex Jacobi operator; this is the restriction of the higher order Jacobi operator to the set of
complex $2$-planes. The following result is well known, for example see
\cite{Gilkey-book}. Conditions (2), (3) of the Lemma simply mean that the operator under
consideration is complex linear. 

\begin{lemma}\label{lem-1.1}
We say that $R$ and  $J$ are {\rm compatible}
if any of the following equivalent conditions are satisfied:
\begin{enumerate}
\item $R(x,y,z,t)=R(Jx,Jy,Jz,Jt)$ for all $x,y,z,t\in V$.
%\item $\mathcal{J}_R(x)J=J\mathcal{J}_R(x)$ for all $x\in V$. (I think this condition is stronger)
\item $\mathcal{J}_R(\pi_x)J=J\mathcal{J}_R(\pi_x)$ for all $x\in
S(V,\langle\cdot,\cdot\rangle)$. \item
$\mathcal{R}(x,Jx)J=J\mathcal{R}(x,Jx)$ for all $x\in
S(V,\g\cdot\cdot)$.
\end{enumerate}\end{lemma}

Note that the curvature $R$ and the almost complex structure $J$ of a K\"ahler manifold are compatible. Thus
this is a very natural geometric condition.

\subsection{The complex Osserman condition}\label{sect-1.5}
Instead of the $2$-Osserman condition where the eigenvalues are
constant on the Grassmannian of $2$-planes, we consider a natural
weaker condition with constant eigenvalues on the space of
holomorphic planes, $\mathbb{CP}(V,\langle\cdot,\cdot\rangle,J)$.

\begin{de}\label{defn-1.2}
Let $\mathcal{V}:=(V,\langle\cdot,\cdot\rangle,J,R)$. We say that $\mathcal{V}$ is a {\rm complex model} if
$\langle\cdot,\cdot\rangle$ is a positive definite inner product on $V$, if
$J$ is an Hermitian almost complex structure on $(V,\g{\cdot}{\cdot})$, and if $R$ is an algebraic
curvature tensor on $(V,\langle\cdot,\cdot\rangle)$. We say that $\mathcal{V}$ is {\rm complex Osserman} if
\begin{enumerate}
\item $\mathcal{V}$ is a complex model.
\item $J$ and $R$ are compatible, i.e. 
$\mathcal{J}_R(\pi_x)$ is complex linear for all $x\in S(V,\langle\cdot,\cdot\rangle)$.
\item The eigenvalues of $\mathcal{J}_R(\pi_x)$ are constant on
$\mathbb{CP}(V,\langle\cdot,\cdot\rangle,J)$.\end{enumerate} We shall also sometimes simply say that $R$ is
{\rm complex Osserman} in this situation.
\end{de}

\subsection{The canonical curvature tensor}\label{sect-1.6}
In addition to the tensor of constant sectional curvature $+1$ defined in
Equation (\ref{eqn-1.b}), it is useful to consider the tensor
\begin{equation}\label{eqn-1.c}
R_\Psi(x,y,z,t):=\g{x}{\Psi t}\g{y}{\Psi z}-\g{x}{\Psi z}\g{y}{\Psi
t}-2\g{x}{\Psi y}\g{z}{\Psi t}
\end{equation}
where $\Psi$ is a skew-symmetric endomorphism of
$(V,\langle\cdot,\cdot\rangle)$. Such tensors play an important
role in studying the space of all algebraic curvature tensors. For
example, Fiedler \cite{Fie} has shown that tensors of this form
span the space of all algebraic curvature tensors. In this paper,
we shall study tensors of this form where the endomorphism in
question defines an Hermitian almost complex structure on
$(V,\langle\cdot,\cdot\rangle)$. We note for future reference that
\begin{equation}\label{eqn-1.d}
\mathcal{J}_{R_0}(x)y=y-\langle y,x\rangle\quad\text{and}\quad
\mathcal{J}_{R_\Psi}(x)y=3\langle y,\Psi x\rangle\Psi x\,.
\end{equation}

\subsection{Algebraic curvature tensors given by Clifford
families}\label{sect-1.7} We say that a set
$\mathcal{F}=\{J_1,\dots,J_\kappa\}$ of Hermitian almost complex
structures on $(V,\langle\cdot,\cdot\rangle)$ is a {\it Clifford
family of rank $\kappa$} if they are subject to the commutation
rules
\[
J_iJ_j+J_jJ_i=-\delta_{ij}\operatorname{id}.
\]
We say that a model $(V,\langle\cdot,\cdot\rangle,R)$ is {\it given by a
Clifford family $\mathcal{F}$ of rank $\kappa$} if there exist constants $c_i$
with $c_i\ne0$ for $1\le i\le\kappa$ so that
\begin{equation}\label{eqn-1.e}
R=c_0R_0+c_1R_{J_1}+...+c_\kappa R_{J_\kappa}\,.
\end{equation}
We shall also sometimes say that $R$ is given by a Clifford family in this setting.
The relations of Equation (\ref{eqn-1.d}) yield that:
\begin{equation}\label{eqn-1.f}
\mathcal{J}_R(x)y=c_0\{y-\langle y,x\rangle x\}
+3c_1\langle y,J_1x\rangle J_1x+...+3c_\kappa\langle y,J_\kappa x\rangle J_\kappa x\,.
\end{equation}
From this it follows immediately that
\begin{equation}\label{eqn-1.g}
\begin{array}{l}
\ja(\pi_x)y=c_0\{2y-\g{y}{x}x-\g{y}{Jx}Jx\}\\
\phantom{\ja(\pi_x)y}+\displaystyle\sum_{i=1}^\kappa
3c_i\{\g{y}{J_ix}J_ix+\g{y}{J_iJx}J_iJx\}\,.
\vphantom{\vrule height 11pt}
\end{array}
\end{equation}

\subsection{Reparametrizing Clifford families}\label{sect-1.8}
Let $A=(A_{ij})\in O(\kappa)$ be an orthogonal
matrix. Set
$$\mathcal{\tilde F}:=\{\tilde
J_i=A_{i1}J_1+\dots+A_{i\kappa}J_\kappa\}\,.$$
This new Clifford family is said to be a {\it
reparametrization} of $\mathcal{F}$; this defines an
equivalence relation on the collection of Clifford families.

\subsection{Summary of results}\label{sect-1.9}
In this paper we begin the study of complex Osserman manifolds by
concentrating on the analysis of complex Osserman models. In Section \ref{preliminaries}, we give
necessary and sufficient conditions so that a model $\mathcal{V}=(V,\langle\cdot,\cdot\rangle,J,R)$ is
complex Osserman and we show that $R$ is necessarily Einstein if $\mathcal{V}$ is complex Osserman. We also
give a topological result in Theorem
\ref{thm-2.4} which controls the eigenvalue structure of $\mathcal{J}_R(\pi_x)$ if $\mathcal{V}$ is complex
Osserman. In Section
\ref{sect-3}, we recall results of Adams on the existence of Clifford families and discuss some
reparametrization results. We also present some examples of
complex Osserman models and show Theorem
\ref{thm-2.4} is sharp.

Work of Nikolayevsky shows that any Osserman model $(V,\langle\cdot,\cdot\rangle,R)$
is given by a Clifford family except in dimension $16$. 
We divide our study into two cases depending on the rank $\kappa$ 
of the structure in question. 

We study the case $\kappa>3$ in Section \ref{sect-4} and show:

\begin{theorem}\label{thm-1.3}
Let $\mathcal{V}:=(V,\langle\cdot,\cdot\rangle,J,R)$ where
$R=c_0R_0+c_1 R_{J_1}+...+c_\kappa R_{J_\kappa}$ is given by a Clifford family of rank $\kappa\geq 4$ on a
vector space $V$ of dimension $n$. The following assertions hold:
\begin{enumerate}
\item Let $c_0=0$. If $\kappa=4,5$, assume $n\geq 2^\kappa$ and, if
$\kappa\geq 6$, assume $n\geq \kappa(\kappa-1)$. Then
$\mathcal{V}$ is not complex Osserman. \item Let $c_0\neq 0$. If
$\kappa=4$ assume $n\geq 32$, if $\kappa=5,6,7$ assume $n\geq
2^\kappa$, if $\kappa\geq 8$ assume $n\geq \kappa(\kappa-1)$. Then
$\mathcal{V}$ is not complex Osserman.
\end{enumerate}
\end{theorem}
Note that, as a consequence of Lemma~\ref{lem-3.1} below, the
hypothesis $n\geq \kappa(\kappa -1)$ in Theorem \ref{thm-1.3} is
not a restriction when $\kappa \geq 16$. Consequently, 
there are only a finite number of possibly exceptional dimensions and ranks when $\kappa\ge4$.

Section \ref{sect-5} is devoted to the study of Clifford families
of lower rank. Results in this section are summarized in the
following theorem:
\begin{theorem}\label{thm-1.4}
Let $\mathcal{V}:=(V,\langle\cdot,\cdot\rangle,J,R)$. Let
$\mathcal{F}=\{J_i\}$ be a Clifford family on a vector space $V$ of dimension
 $n$. Let $c_i\ne0$ be given for $1\le i\le\kappa$ where $\kappa\le 3$.
\begin{enumerate}
\item \underline{Rank $\kappa=0$}.
Let $R=c_0 R_0$. Then $\mathcal{V}$ is complex Osserman. \item
\underline{Rank $\kappa=1$}. Let $R=c_0R_0+c_1R_{J_1}$.
\begin{enumerate}
\item[(a)] If $c_0=0$, then $\mathcal{V}$ is complex Osserman if
and only if $J J_1=\pm J_1 J$. \item[(b)] If $c_0\ne0$, then $\mathcal{V}$ is
complex Osserman if and only if $J=\pm J_1$ or $JJ_1=-J_1J$.
\end{enumerate}
\item \underline{Rank $\kappa=2$}.
Let $R=c_0R_0+c_1R_{J_1}+c_2R_{J_2}$. Then $\mathcal{V}$ is complex Osserman if and
only if there exists a reparametrization $\{\tilde J_1,\tilde J_2\}$ of
$\mathcal{F}$ so that one has $R=c_0R_0+\tilde c_1 R_{\tilde J_1}+\tilde c_2
R_{\tilde J_2}$ and so that one of the following holds:
\begin{enumerate}
\item[(a)] $c_0=0$, $J\tilde J_1=\tilde J_1 J$ and $J\tilde
J_2=-\tilde J_2 J$. \item[(b)] Either $J=\tilde J_1$ or $J=\tilde
J_1 \tilde J_2$.
\end{enumerate}
\item \underline{Rank $\kappa=3$}. Let $R=c_0R_0+c_1R_{J_1}+c_2R_{J_2}+c_3R_{J_3}$.
\begin{enumerate}
\item[(a)] Assume $n\geq 12$. If $c_0=0$, then $\mathcal{V}$ is
complex Osserman if and only if there exists a reparametrization
$\{\tilde J_1,\tilde J_2,\tilde J_3\}$ of $\mathcal{F}$ so that one
has $R=\tilde c_1R_{\tilde J_1}+\tilde c_2 R_{\tilde J_2}+\tilde c_3 R_{\tilde J_3}$
and that $J=\tilde J_1$ or $J=\tilde J_2\tilde J_3$.
\item[(b)] Assume
$n\geq 16$. If $c_0\ne0$, then $\mathcal{V}$ is complex
Osserman if and only if there exists a reparametrization
$\{\tilde J_1,\tilde J_2,\tilde J_3\}$ of $\mathcal{F}$ so that one has
$R=c_0R_0+\tilde c_1R_{\tilde J_1}+\tilde c_2R_{\tilde J_2}+\tilde c_3R_{\tilde J_3}$,
$J=\tilde J_1$, and
$\tilde J_1\tilde J_2\tilde J_3=\operatorname{id}$.
\end{enumerate}
\end{enumerate}
\end{theorem}

\begin{remark}\label{rmk-1.5}\rm
From Theorem \ref{thm-1.4} we obtain the following geometric
conclusions:
\begin{enumerate} \item Let $(M,g)$ be a manifold of
constant sectional curvature. Then $(M,g)$ is complex Osserman
with respect to any Hermitian almost complex structure $J$.

\item Let $(M,g,J)$ be a K{\"a}hler manifold which has constant holomorphic
sectional curvature. Then $(M,g,J)$ is complex Osserman with
respect to $J$.

\item Let $(M,g,\{J_1,J_2,J_3\})$ be a quaternionic K{\"a}hler
manifold which has constant quaternionic sectional curvature, where
$\{J_1,J_2,J_3\}$ forms a locally defined quaternionic structure.
Then, for any $J\in \operatorname{Span}\{J_1,J_2,J_3\}$, $(M,g,J)$ is complex
Osserman.
\end{enumerate}
Note that if $(M,g)$ is Osserman of dimension different from $16$,
then it is isometric to one of these three examples or is flat \cite{Chi,nikolayevsky-1,nikolayevsky-2}.
\end{remark}

\section{Algebraic preliminaries}\label{preliminaries}

In this section we present some foundational results. Our first
result is the well known observation:
\begin{lemma}\label{lem-2.1}
Let $\mathcal{V}_i:=(V,\langle\cdot,\cdot\rangle,R_i)$, with
$i=1,2$, be models. If $\mathcal{J}_{R_1}(x)=\mathcal{J}_{R_2}(x)$
for all $x$ in $V$, then $R_1=R_2$.
\end{lemma}

What is perhaps somewhat surprising is that this observation fails
for the complex Jacobi operator as we shall see in
Theorem \ref{thm-3.6}. Let $\operatorname{Spec}\{\ja(\pi_x)\}$
be the spectrum of $\ja(\pi_x)$ and let $E_\lambda(\pi_x)$ be the
eigenspace associated to the eigenvalue $\lambda$ of $\ja(\pi_x)$.
Since $\mathcal{J}_R(\pi_x)$ is self-adjoint,
$\mathcal{J}_R(\pi_x)$ is diagonalizable with real eigenvalues.
Thus we have an orthogonal direct sum decomposition
$$V=\oplus_\lambda E_\lambda(\pi_x)$$
for any $x\in S(V,\langle\cdot,\cdot\rangle)$. The following lemma is an immediate
consequence of Lemma~\ref{lem-1.1} and provides a criterion for complex Osserman
curvature tensors:
\begin{lemma}\label{lem-2.2}
$\mathcal{V}=(V,\g{\cdot}{\cdot},J,R)$ is complex Osserman if and only if
\begin{enumerate}
\item $J E_\lambda(\pi_x)=E_\lambda(\pi_x)$ for all $\pi_x\in
\mathbb{CP}(V,\langle\cdot,\cdot\rangle,J)$ and $\lambda\in
\operatorname{Spec}\{\ja(\pi_x)\}$. \item
$\operatorname{Spec}\{\ja(\pi_x)\}=\operatorname{Spec}\{\ja(\pi_y)\}$ for all $\pi_x,\pi_y\in
\mathbb{CP}(V,\langle\cdot,\cdot\rangle,J)$.
\end{enumerate}
\end{lemma}

An model $(V,\langle\cdot,\cdot\rangle,R)$ is said to be Einstein if
$\rho(\cdot,\cdot)=c\,\g\cdot\cdot$ for a constant $c$, where by
$\rho$ we denote the Ricci tensor. In general, $p$-Osserman
models are Einstein. This result generalizes to become:

\begin{lemma}\label{lem-2.3}
Let $\mathcal{V}=(V,\langle \cdot,\cdot \rangle,J,R)$ be complex
Osserman. Then $\mathcal{V}$ is Einstein.
\end{lemma}
\begin{proof}
Assume that $\mathcal{V}$ is complex Osserman. Let $x\in
S(V,\langle\cdot,\cdot\rangle)$. As $R$ is compatible,
$R(y,Jx,Jx,z)=R(Jy,x,x,Jz)$ and thus
$\mathcal{J}_R(Jx)=-J\mathcal{J}_R(x)J$. Consequently
$$\rho(x,x)=\operatorname{Tr}\{\mathcal{J}_R(x)\}=\operatorname{Tr}\{\mathcal{J}_R(Jx)\}
=\textstyle\frac12\operatorname{Tr}\{\mathcal{J}_R(\pi_x)\}=\textstyle\frac12
\sum_\lambda\lambda\dim\{E_\lambda(\pi_x)\}$$ is independent of
$x\in S(V,\langle\cdot,\cdot\rangle)$. This implies
$\rho(\cdot,\cdot)=c\,\g\cdot\cdot$. Consequently $\mathcal{V}$ is
Einstein.
\end{proof}

Methods of algebraic topology can be used to control the
eigenvalue structure of a complex Osserman model. In particular, no
more than $3$ distinct eigenvalues may occur.

\begin{theorem}\label{thm-2.4}
Let $\mathcal{V}=(V,\langle \cdot,\cdot \rangle,J,R)$ be complex
Osserman. If $\jp$ is not a multiple of the identity (i.e. if $\mathcal{J}_R(\pi_x)$ has at least 2
distinct eigenvalues), then:
\begin{enumerate}
\item If $n\equiv 2\; (mod\; 4)$, there are $2$ eigenvalues with
multiplicities $(n-2,2)$. \item If $n\equiv 0\; (mod\; 4)$, then
one of the following holds:
\begin{enumerate}
\item There are $2$ eigenvalues with multiplicities $(n-2,2)$.
\item There are $2$ eigenvalues with multiplicities $(n-4,4)$.
\item There are $3$ eigenvalues with multiplicities $(n-4,2,2)$.
\end{enumerate}
\end{enumerate}
\end{theorem}
\begin{proof} Let $\mathbb{V}:=\mathbb{CP}(V,\langle\cdot,\cdot\rangle,J)\times V$ be the trivial bundle
over projective space.
Lemma \ref{lem-2.2} shows that the
eigenspaces
\[
E_{\lambda_i}(\pi):=\{v\in V: \ja(\pi)v=\lambda_i v\}
\]
have constant rank and patch together to define smooth vector
bundles $E_{\lambda_i}(\pi)$ over $\mathbb{CP}(V,\langle\cdot,\cdot\rangle,J)$ where
$\{\lambda_0,...,\lambda_k\}$ denote the distinct eigenvalues of $\mathcal{J}_R(\pi)$ for any, and hence
for all,
$\pi\in\mathbb{CP}(V,\langle\cdot,\cdot\rangle,J)$. This gives the following direct sum decomposition
\[
\mathbb{V}=E_{\lambda_0}\oplus\dots\oplus E_{\lambda_k}.
\]
This decomposition is in the category of complex vector bundles
since the eigenbundles are invariant under $J$.

A sub-bundle $E$ of $\mathbb{V}$ is said to be a {\it geometrically symmetric
vector bundle} if for all complex lines $\sigma,\tau$ in
$\mathbb{CP}(V,\langle\cdot,\cdot\rangle,J)$,
$\tau
\subset E(\sigma)$ implies $\sigma \subset E(\tau)$.
Let $\lambda_{\operatorname{min}}$ be the minimal eigenvalue of $\mathcal{J}_R(\pi_x)$. We
then have the following chain of equivalences for unit vectors $x$ and $y$:
\begin{eqnarray*}
&&y\in E_{\lambda_{\operatorname{min}}}(\pi_x)\\
&\Leftrightarrow&R(y,x,x,y)+R(y,Jx,Jx,y)=\lambda\\
&\Leftrightarrow&R(x,y,y,x)+R(Jx,y,y,Jx)=\lambda\\
&\Leftrightarrow&
R(x,y,y,x)+R(x,Jy,Jy,x)=\lambda\\
&\Leftrightarrow&x\in E_{\lambda_{\operatorname{min}}}(\pi_y)\,.
\end{eqnarray*}
This implies that the bundle $E_{\lambda_{\operatorname{min}}}$ is
geometrically symmetric. The desired result now follows from
results in \cite{G01} concerning geometrically symmetric bundles;
these results generalize earlier results of Glover et al.
\cite{GHS}.\end{proof}

We shall show that this result is sharp in Remark \ref{rmk-3.5}
below by showing that all the possibilities can be realized.

\section{Clifford families and associated curvature tensors}\label{sect-3}

A Clifford family $\mathcal{F}=\{J_1,J_2,J_3\}$ of rank $3$ is
called a {\it quaternion structure} if $J_1 J_2=J_3$. Note that $V$
admits a quaternion structure if and only if
$\dim\{V\}$ is divisible by $4$. One defines the {\it Adams number} $\nu(n)$ by setting
$\nu(1)=0$,
$\nu(2)=1$, $\nu(4)=3$, $\nu(8)=7$, $\nu(16 r)=\nu(r)+8$ and
$\nu(m 2^s)=\nu(2^s)$ for $m$ odd. One then has the following well known result of Atiyah
et al. \cite{ABS} which is closely related to work of Adams \cite{Adams} concerning vector
fields on spheres:

\begin{lemma}\label{lem-3.1}
There exists a Clifford family of rank $\kappa$ on $V$ if and only if $\kappa\le\nu(n)$.
\end{lemma}
We now present a useful technical result:

\begin{lemma}\label{lem-3.2}
Let $V$ and $W$ be vector spaces and let
$\mathcal{T}=\{T_1,\dots,T_\kappa\}$ be a family of linear
maps $T_i:V\rightarrow W$. Assume there is an integer $\mu$ so that
for any set of constants $a_i$, not all of which are zero, one has
$\operatorname{Rank}\{a_1T_1+\dots+a_\kappa T_\kappa\}\geq \mu$. Then the following assertions hold:
\begin{enumerate}
\item If $\mu\geq \kappa$, there exists $x\in V$ so that
$\{T_1x,\dots,T_\kappa x\}$ is a set of linearly independent
vectors. 
\item If $\mu\geq 2\kappa$, there exists $x,y\in V$ so
that $\{T_1x,\dots,T_\kappa x,T_1y,\dots,T_\kappa y\}$ is a set of
linearly independent vectors. 
\item Let $T:V\longrightarrow W$ be
a linear map so that $Tx\in \operatorname{Span}\{T_1x,\dots,T_\kappa x\}$ for
all $x\in V$. If $\mu\ge2\kappa$, then $T\in \operatorname{Span}\{T_1,\dots,T_\kappa\}$.
\end{enumerate}
\end{lemma}
\begin{proof} In order to prove Assertion (1), suppose $\mu\geq
\kappa$. For a given $x\in V$, choose $r(x)$ maximal so that
$\{T_1x,\dots,T_rx\}$ is a linearly independent set of $r$ vectors.
Take $x\in V$ so that $r(x)$ is maximal. If $r(x)=\kappa$, then clearly
Assertion (1) holds. Suppose $r(x)<\kappa$. We argue for a
contradiction. Choose $(a_1,\dots,a_r)$ so that
$a_1T_1x+\dots+a_rT_rx+T_{r+1}x=0$ and let
$$S:=a_1T_1+\dots+a_rT_r+T_{r+1}\,.$$
As $\operatorname{Rank}\{S\}\geq
\mu\geq\kappa$, there is $y\in V$ so that $\{T_1x,\dots,T_rx,Sy\}$ is a set of
$r+1$ linearly independent vectors. Hence, by continuity, there
exists $\epsilon>0$ such that
$\{T_1(x+\epsilon y),\dots,T_r(x+\epsilon y),Sy\}$ 
is a set of $r+1$ linearly independent
vectors. Consequently 
$\{T_1(x+\epsilon y),\dots,T_r(x+\epsilon y),T_{r+1}(x+\epsilon y)\}$ 
also is a set of $r+1$ linearly independent vectors. 
Therefore $r(x+\epsilon y)\geq r+1$ which contradicts the choice of $x$. This
contradiction establishes Assertion (1).

Now suppose that $\mu\geq 2\kappa$. By Assertion (1) we may choose
$x\in V$ so that $\{T_1x,\dots,T_\kappa x\}$ is a linearly
independent set of $\kappa$ vectors. Consider the vector space
$W_0:=\operatorname{Span}\{T_1x,\dots,T_\kappa x\}$ and let
$\pi:W\longrightarrow W/W_0$ be the natural projection. We apply Assertion (1) to
the linear maps $\bar T_i:=\pi T_i:V\longrightarrow W/W_0$ with $\bar\mu=\mu-\kappa\ge\kappa$
to complete the proof of Assertion (2).

We complete the proof by establishing Assertion (3). By assumption, for every 
$z\in V$, there exist coefficients $a_i(z)$ so that
$Tz=a_1(z)T_1z+...+a_\kappa(z)T_\kappa z\,.$
To show that $T\in \operatorname{Span}\{T_1,\dots,T_\kappa\}$, we must show that
the coefficients can be chosen to be independent of $z$.

By Assertion (2), there are vectors $x,y \in S(V,\langle\cdot,\cdot\rangle)$
so $\{T_1x,...,T_\kappa x,T_1y,...,T_\kappa y\}$ is a collection of $2\kappa$ linearly
independent vectors. Then, by continuity, this remains true on some
open neighborhoods $\mathcal{O}_x$ and $\mathcal{O}_y$ of $x$ and
$y$, respectively. Let $z\in\mathcal{O}_x$ and let
$t\in\mathcal{O}_y$. We may then express:
\begin{eqnarray*}
T(z+t)&=&\sum_{i=1}^\kappa a_i(z+t)T_i(z+t)=\sum_{i=1}^\kappa a_i(z+t)(T_iz+T_it)\\
      &=&Tz+Tt=\sum_{i=1}^\kappa\{a_i(z)T_iz+a_i(t)T_it\}\,.
\end{eqnarray*}
Since the vectors $\{T_1z,...,T_\kappa z,T_1t,...,T_\kappa z\}$
are linearly independent, this implies $a_i(z)=a_i(z+t)=a_i(t)$
for $z\in \mathcal{O}_x$ and $t\in \mathcal{O}_y$. Thus, for
$a_i:=a_i(t)$,
$$Tz=\sum_{i=1}^\kappa
a_iT_iz\quad\text{for all}\quad z\in\mathcal{O}_x\,.$$ This
polynomial identity holds on a non-empty open set and thus holds
on all $V$. This establishes Assertion (3).
\end{proof}

We specialize this result for Clifford families.

\begin{corollary}\rm \label{cor-3.3}
Let $\mathcal{F}:=\{J_1,\dots,J_\kappa\}$ be a Clifford family of rank $\kappa$
on a vector space of dimension $n$.
\begin{enumerate}
\item Suppose that $n\ge\kappa$. Then there exists $x$ in $V$ so that the set $\{J_ix\}_{1\le i\le\kappa}$
consists of $\kappa$ linearly independent vectors.
\item Suppose that $n\ge2\kappa$. Then there exist $x$ and $y$ in $V$ so that the set $\{J_ix,J_iy\}_{1\le i\le\kappa}$ consists of $2\kappa$ linearly independent vectors. Furthermore, if
$Tx\in\operatorname{Span}_{1\le i\le\kappa}\{J_ix\}$ for all $x$ in $V$, then
$T\in\operatorname{Span}_{1\le i\le\kappa}\{J_i\}$.
\item Suppose that $n\ge\kappa(\kappa-1)$. Then
there exists $x$ in $V$ so that the set
$\{J_jJ_kx\}_{1\le j<k\le\kappa}$ consists of
$\frac12\kappa(\kappa-1)$ linearly independent vectors.
\item Suppose that $n\ge2\kappa(\kappa-1)$. Then
there exist $x$ and $y$ in $V$ so that the set
$\{J_jJ_kx,J_jJ_ky\}_{1\le j<k\le\kappa}$ consists of
$\kappa(\kappa-1)$ linearly independent vectors. Furthermore, 
if $Tx\in\operatorname{Span}_{1\le j<k\le\kappa}\{J_jJ_kx\}$ for all $x$ in $V$, then
$T\in\operatorname{Span}_{1\le j<k\le\kappa}\{J_jJ_k\}$.
\end{enumerate}
\end{corollary}

\begin{proof}
One verifies that $(a_1J_1+...+a_\kappa J_\kappa)^2=-(a_1^2+...+a_\kappa^2)\operatorname{id}$
and thus one has that 
$\operatorname{Rank}(a_1J_1+...+a_\kappa J_\kappa)=n$ if any coefficient is
non-zero. Assertions (1) and (2) now follow from Lemma \ref{lem-3.2}. If not all the coefficients vanish,
one shows similarly that:
$$
\operatorname{Rank}\left(\sum_{j=1}^{\kappa-1}\sum_{k=j+1}^\kappa 
     a_{jk}J_jJ_k\right)\ge\frac n2\,.
$$
The remaining assertions of the Lemma now follow.\end{proof}

We now describe some general properties of models given by Clifford families. We adopt the notation of
Equations (\ref{eqn-1.b}) and (\ref{eqn-1.c}).

\begin{lemma}\label{lem-3.4}
\ \begin{enumerate}
\item  Suppose that $J$ is an Hermitian almost complex structure on $(V,\langle\cdot,\cdot\rangle)$.
Then $\mathcal{V}:=(V,\langle\cdot,\cdot\rangle,J,c_0R_0+c_1R_J)$ is complex Osserman.
\item Suppose that $\{J_1,J_2,J_3\}$ is an Hermitian quaternion structure on
$(V,\langle\cdot,\cdot\rangle)$. Then
$\mathcal{V}:=(V,\langle\cdot,\cdot\rangle,J_1,c_0R_0+c_1R_{J_1}+c_2R_{J_2}+c_3R_{J_3})$ is complex
Osserman.
\item Let $\f:=\{J_1,\dots,J_\kappa\}$ be a Clifford family and let
$\tilde\f:=\{\tilde J_1,\dots,\tilde J_\kappa\}$ be a
reparametrization of $\f$. Then
$R_{J_1}+\dots+R_{J_\kappa}=R_{\tilde J_1}+\dots+R_{\tilde
J_\kappa}$.\end{enumerate}
\end{lemma}
\begin{proof}
Let $\mathcal{V}$ be as in Assertion (1). We use Equation (\ref{eqn-1.g}) to see that:
$$
\mathcal{J}_{R}(\pi_x)y=\left\{\begin{array} {rl}
(c_0+3c_1)y\!&\! \mbox{if}\; y\in \operatorname{Span}\{x,Jx\},\\
2c_0\!&\! \mbox{if}\; y\bot\, \operatorname{Span}\{x,Jx\}.
\end{array}\right.
$$
Hence $J$ and $\jp$ commute and the eigenvalues are constant. Thus $\mathcal{V}$ is complex Osserman 
by Lemma \ref{lem-2.2}; the proof of Assertion (2) is similar and follows from a
calculation in this instance that:
$$
\mathcal{J}_{R_J}(\pi_x)y=\left\{\begin{array} {rl}
(c_0+3c_1)y\!&\! \mbox{if}\; y\in \operatorname{Span}\{x,J_1x\},\\
(2c_0+3c_2+3c_3)y\!&\! \mbox{if}\; y\in \operatorname{Span}\{J_2x,J_3x\},\\
2c_0\!&\! \mbox{if}\; y\bot\, \operatorname{Span}\{x,J_1x,J_2x,J_3x\}.
\end{array}\right.$$

We complete the proof by verifying that Assertion (3) holds.
If $x\in S(V)$, then the vectors $\{J_1x,...,J_\kappa x\}$ form an
orthonormal set. Let $\sigma_{\mathcal{F}}x$ be orthogonal
projection on the subspace
$$
  S_1^{\mathcal{F}}(x):=\operatorname{Span}\{J_1x,...,J_\kappa x\}\,.
$$
We then have $\sum_i\langle x,J_ix\rangle
J_ix=\sigma_{\mathcal{F}}x$.
 Let
$R=R_{J_1}+...+R_{J_\kappa}$. By Equation
(\ref{eqn-1.f}),
$\mathcal{J}_R(x)=3\sigma_{\mathcal{F}}(x)$.
If $\tilde{\mathcal{F}}$ is a reparametrization of $\mathcal{F}$,
then $S_1^{\mathcal{F}}(x)=S_1^{\tilde{\mathcal{F}}}(x)$.
Consequently $\mathcal{J}_R(x)=\mathcal{J}_{\tilde R}(x)$ so by Lemma \ref{lem-2.1}, $R=\tilde R$.
\end{proof}

\begin{remark}\label{rmk-3.5}\rm
Theorem \ref{thm-2.4} places restrictions on the possible eigenvalue multiplicities
of the complex Jacobi operator defined by a complex Osserman model. We may use Lemma 
\ref{lem-3.4} to show that in fact all these possibilities occur.
Suppose first that the dimension $n$ of $V$ is even. Let $J$ be an Hermitian almost complex
structure on $(V,\langle\cdot,\cdot\rangle)$.
\begin{enumerate}
\item If $R=3R_0+R_J$, then $J_R(\pi_x)=6\operatorname{id}$.
\item If $R=R_0+R_J$, then the eigenvalues of $J_R(\pi_x)$ are $(2,4)$ and the eigenvalue multiplicities
are
$(n-2,2)$.
\end{enumerate}
If $n$ is divisible by $4$, there are additional eigenvalue multiplicities which
can be realized. Let $\{J_1,J_2,J_3\}$ be a quaternion structure on 
$(V,\langle\cdot,\cdot\rangle)$ and let $J=J_1$.
\begin{enumerate}
\item If $R=3R_0+3R_{J_1}+R_{J_2}+R_{J_3}$, then $\mathcal{J}_R(\pi_x)$ are $(6,12)$ and
the eigenvalue multiplicities are
$(n-4,4)$.
\item If $R=R_0+R_{J_1}+R_{J_2}+R_{J_3}$, then the eigenvalues of $\mathcal{J}_R(\pi_x)$ are $(2,4,8)$ and
the eigenvalue multiplicities are
$(n-4,2,2)$.
\end{enumerate}\end{remark}

Lemma \ref{lem-2.1} shows that the Jacobi operator determines the full curvature
tensor, i.e. that if $\mathcal{J}_R(x)=0$ for all $x$ in $V$, then $R=0$. Similarly,
the higher order Jacobi operator determines the full curvature operator. To see that
this is true, one may argue as follows. Let $2\le p\le n-1$. Assume that $\mathcal{J}_R(\sigma)=0$
for every $p$-plane $\sigma$. Let $x$ and $y$ be unit vectors in $V$. Choose additional
unit vectors $\{e_2,...,e_p\}$ so that $\{x,e_2,...,e_p\}$ is an orthonormal basis
for a $p$-plane $\sigma_x$ and so that $\{y,e_2,...,e_p\}$ is an orthonormal basis
for a $p$-plane $\sigma_y$. Then
$$
\mathcal{J}_R(x)y=(\mathcal{J}_R(x)-\mathcal{J}_R(y))y
=(\mathcal{J}_R(\sigma_x)-\mathcal{J}_R(\sigma_y))y=0\,.
$$
This shows that $\mathcal{J}_R=0$ and hence $R=0$ by Lemma \ref{lem-2.1}.

However an analogous property does not hold for the complex Jacobi operator. This is, 
perhaps, to be expected on dimensional grounds. The domain of the usual Jacobi 
operator is $V$ which is $n$-dimensional. The domain of the higher order Jacobi operator
is the dimension of the $p$-dimensional Grassmannian which has dimension 
greater than $n$
for $2\le p\le n-2$. However, the domain of the complex Jacobi operator is
$\mathbb{CP}(V,\langle\cdot,\cdot\rangle,J)$ which is $n-2$ dimensional. One has the following result:

\begin{theorem}\label{thm-3.6}
Let $V$ be a vector space of dimension $n$. Assume $n$ is divisible by $4$ 
and that $n$ is at least $8$.
Then there exists
a model $\mathcal{V}=(V,\langle\cdot,\cdot\rangle,J,R)$ which is
complex Osserman, which is not Osserman, which is not given by
a Clifford family, and which has $\mathcal{J}_R(\pi_x)=0$ for all
$x$.\end{theorem}

\begin{proof} Since the dimension of $V$ is divisible by $4$, we can
find a quaternion structure $\{K_1,K_2,K_3\}$ on $V$. Since $n\ge8$,
we may take a non-trivial decomposition of $V$ as a quaternion module
in the form $V=V_+\oplus V_-$. Define a new Clifford family on $V$ which is not
a quaternion structure by setting $J_1:=K_1$, $J_2:=K_2$, and $J_3:=\mp J_1J_2$
on $V_\pm$. We then have $J_1J_2J_3x=\pm x$ for $x\in V_\pm$.
Define
$$R:=R_{J_2}-R_{J_1J_2}-R_{J_3}+R_{J_1J_3}\,.$$
Let $x_\pm\in S(V_\pm)$. Equation (\ref{eqn-1.d}) yields that:
$$
\mathcal{J}_R(x_+)y=\left\{\begin{array}{rll}
6y&\text{ if }&y\in\operatorname{Span}\{J_2x_+\}=\operatorname{Span}\{J_1J_3x_+\},\\
-6y&\text{ if }&y\in\operatorname{Span}\{J_3x_+\}=\operatorname{Span}\{J_1J_2x_+\}\,.
\end{array}\right.
$$
On the other hand, if we take $x_0=(x_++x_-)/\sqrt{2}$, then
$$
\mathcal{J}_R(x_0)y=\left\{\begin{array}{rll}
3y&\text{ if }&y\in\operatorname{Span}\{J_2x_0,J_1J_3x_0\}=\operatorname{Span}\{J_2x_+,J_2x_-\},\\
-3y&\text{ if }&y\in\operatorname{Span}\{J_1J_2x_0,J_3x_0\}=\operatorname{Span}\{J_3x_+,J_3x_-\}\,.
\end{array}\right.$$
This shows that $\mathcal{V}$ is not Osserman. As any model given by a Clifford family is necessarily
Osserman, $\mathcal{V}$ is not given by a Clifford family. On the
other hand, the complex Jacobi operator with respect to $J=J_1$ is
given by
\begin{eqnarray*}
\mathcal{J}_R(\pi_x)y&=&3\langle y,J_2x\rangle J_2x+3\langle y,J_2J_1x\rangle J_2J_1x
-3\langle y,J_1J_2x\rangle J_1J_2x\\
&-&3\langle y,J_1J_2J_1x\rangle J_1J_2J_1x
 -3\langle y,J_3x\rangle J_3x-3\langle y,J_3J_1x\rangle J_3J_1x\\
&+&3\langle y,J_1J_3x\rangle J_1J_3x+3\langle y, J_1J_3J_1x\rangle J_1J_3J_1x\\
&=&0\,.
\end{eqnarray*}
This shows $\mathcal{J}_R(\pi_x)=0$ for all $x$ as desired. Thus $\mathcal{V}$ is complex
Osserman.\end{proof}

\section{Curvature and higher order Clifford
families}\label{sect-4}

In this section, we establish Theorem \ref{thm-1.3} by studying
models with
$$R=c_0R_0+c_1R_{J_1}+...+c_\kappa R_{J_\kappa}$$
where $\{J_1,...,J_\kappa\}$ is a Clifford family of rank $\kappa\geq 4$ on
$(V,\langle\cdot,\cdot\rangle)$. We remark that the work of \cite{Chi,nikolayevsky-1,nikolayevsky-2} shows
 tensors of this kind do not arise in the geometric context. In  Section \ref{sect-4.1} we study the
case
$c_0=0$ and in Section \ref{sect-4.2} we study the case $c_0\ne0$. We shall always assume that the
constants $c_1$, ..., $c_\kappa$ are  non-zero.

\subsection{Curvature given by a Clifford family with $c_0=0$.}
\label{sect-4.1} Throughout this section we shall assume that
$$R=c_1R_{J_1}+...+c_\kappa R_{J_\kappa}$$
where $c_1$, ..., $c_k$ are non-zero constants and where
$\{J_1,...,J_\kappa\}$ is a Clifford family of rank $\kappa$
on a vector space $V$ of dimension $n$. Let $\mathcal{V}:=(V,\langle\cdot,\cdot\rangle,J,R)$.
We suppose that $\mathcal{V}$ is complex Osserman. We first show that this implies
that $J$ has the form $J=\sum_{i<j}c_{ij}J_iJ_j$. We then derive a contradiction
by studying the eigenvalue structure and by studying the 
coefficients $c_{ij}$. The eigenvalue multiplicity estimates of
Theorem \ref{thm-2.4} will play a crucial role in our analysis.

We shall have to impose certain conditions on $n$; these
conditions are automatic for $\kappa$ large. We begin with a technical result:

\begin{lemma}\label{lem-4.1}
Let $\mathcal{V}:=(V,\g\cdot\cdot,J,R=c_1R_{J_1}+\dots+c_\kappa
R_{J_\kappa})$ be a complex Osserman model where
$\{J_1,...,J_\kappa\}$ is a Clifford family of rank $\kappa$ on a vector space
of dimension $n$. Assume that $\kappa\ge4$ and that $n\ge2\kappa+5$.
If $x\in S(V,\langle\cdot,\cdot\rangle)$, then
\begin{enumerate}
\item $\operatorname{Rank}\{\jp\}\leq 4$.
\item $Jx\in\operatorname{Span}_{i\leq 4,i\neq j}\{J_iJ_jx\}$.
\end{enumerate}
\end{lemma}
\begin{proof}
Equation (\ref{eqn-1.g}) shows
$\operatorname{Rank}\{\jp\}\leq 2\kappa$. Consequently $0$ is an eigenvalue of
multiplicity at least $n-2\kappa\geq 5$.
Theorem \ref{thm-2.4} then shows that $0$ is an eigenvalue of multiplicity
at least $n-4$. Consequently, as desired, $\operatorname{Rank}\{\jp\}\leq 4$.

The vectors $\{J_1x,...,J_\kappa x\}$ form an orthonormal set for $x\in S(V,\langle\cdot,\cdot\rangle)$.
Let $\alpha_i(x):=\g{J_ix}{J_1Jx}$ be the Fourier coefficients of $J_1Jx$. Let
\begin{eqnarray*}
U(x)&:=&\operatorname{Span}\{J_1x,\dots,J_\kappa x,J_1Jx\},\\
V(x)&:=&\operatorname{Span}\{J_2Jx,\dots,J_\kappa Jx\},\\
W(x)&:=&U(x)+V(x).
\end{eqnarray*}
Note that $\operatorname{Range}\{\jp\}\subset W(x)$. If
$\dim\{U(x)\}\leq \kappa$, then $J_1Jx\in\operatorname{Span}_i\{J_ix\}$. Since $J_1Jx\perp J_1x$,
we have that $Jx\in\operatorname{Span}_{i>1}\{J_1J_ix\}$ and Assertion (2) follows.

Suppose on the other hand that $\dim\{U(x)\}=\kappa +1$ or equivalently that
\begin{equation}\label{eqn-4.a}
\alpha_1^2+\dots+\alpha_\kappa^2<1\,.
\end{equation}
Let $\rho$ be the projection on $W(x)/V(x)$. Then
\begin{eqnarray*}
\rho\jp J_ix&=&\rho\{3c_iJ_ix+3c_1\alpha_iJ_1Jx\},\\
\rho\jp J_1Jx&=&\rho\{3c_1JJ_1x+3c_1\alpha_1J_1x+...+3c_\kappa\alpha_\kappa J_\kappa x\}.
\end{eqnarray*}
Hence $\rho\jp=\rho M$ on $U(x)$, where
\[
M:=3 \left(\begin{array}{ccccc}
                c_1&0&\dots&0&c_1\alpha_1\\
                0&c_2&\dots&0&c_2\alpha_2\\
                \dots&\dots&\dots&\dots&\dots\\
                0&0&\dots&c_\kappa&c_\kappa\alpha_\kappa\\
                c_1\alpha_1&c_1\alpha_2&\dots&c_1\alpha_\kappa&c_1
\end{array}\right)
\]
We compute 
$\det(M)=3^{\kappa+1}c_1^2c_2\dots c_\kappa(1-\alpha_1^2-\dots-\alpha_\kappa^2)$. 
Thus by Equation (\ref{eqn-4.a}), $\det(M)\neq 0$ so $M$
is invertible. Consequently,
\[
\dim\{\rho U(x)\}=\dim\{\rho MU(x)\}=\dim\{\rho\jp U(x)\}\leq
\operatorname{Rank}\{\jp\}\leq 4.
\]
The short exact sequence
\[
0\rightarrow V(x)\rightarrow W(x)\rightarrow W(x)/V(x)=\rho
U(x)\rightarrow 0
\]
shows that $\dim\{W(x)\}=\dim\{V(x)\}+\dim\{\rho U(x)\}\leq
(\kappa-1)+4$. Therefore 
\begin{eqnarray*}
\dim\{\operatorname{Span}_{i\leq4}\{J_ix\}\cap\operatorname{Span}_i\{J_iJx\}\}
  &=&4-(\dim\{W(x)\}-\kappa)\\
  &\geq&4+\kappa-(\kappa+3)>0\,.
\end{eqnarray*}
Hence, there exist non-zero constants $a_i$ and $b_j$ so that
\[
a_1J_1x+a_2J_2x+a_3J_3x+a_4J_4x=b_1J_1Jx+\dots+b_\kappa J_\kappa
Jx.
\]
We multiply by $b_1J_1+...+b_\kappa J_\kappa$ to invert this relation and conclude thereby that $Jx\in
\operatorname{Span}\{x,
\{J_iJ_jx\}_{i\leq 4,i\ne j}\}$. Since $Jx\perp x$, we may conclude as desired that
$Jx\in\operatorname{Span}_{i\leq 4,i\ne j}\{J_iJ_jx\}$.
\end{proof}

We continue our study by reducing to the cases $\kappa=4$ and $\kappa=5$:

\begin{lemma}\label{lem-4.2}
Let $\mathcal{V}:=(V,\g\cdot\cdot,J,R=c_1R_{J_1}+\dots+c_\kappa
R_{J_\kappa})$ be a complex Osserman model where
$\{J_1,...,J_\kappa\}$ is a Clifford family of rank $\kappa$ on a vector space
of dimension $n$. Assume that $\mathcal{V}$ is complex Osserman, that $n\geq
\kappa(\kappa-1)$, and that $\kappa\ge4$. Then $\kappa \leq 5$.
\end{lemma}
\begin{proof}
Suppose $\kappa\geq 6$. By Corollary \ref{cor-3.3} we know that
there exists $x\in V$ such that $\{J_iJ_jx\}_{i<j}$ is a linearly
independent set of $\frac12\kappa(\kappa-1)$ vectors. By Lemma \ref{lem-4.1},
\[
Jx=\sum_{1\leq i\leq 6,i<j} a_{ij}(x)J_iJ_jx.
\]
Moreover, the sum may be restricted to $i\leq 4$ and, since the
coefficients $a_{ij}$ are uniquely determined, we get
$a_{56}(x)=0$. By permuting the
role of the indices we may conclude that all the coefficients vanish. As this is not possible,
$\mathcal{V}$ can not be a complex Osserman model.
\end{proof}

The analysis of the cases $\kappa=4$ and $\kappa=5$ to complete the
proof of Theorem \ref{thm-1.3} (1) is a bit technical. We shall
outline the proof but omit details in the interests of brevity.
We assume $\dim(V)\ge16$ throughout.

\begin{lemma}\label{lem-4.3}
Let $\mathcal{V}:=(V,\g\cdot\cdot,J,R=c_1R_{J_1}+\dots+c_\kappa
R_{J_\kappa})$ where
$\{J_1,...,J_\kappa\}$ is a Clifford family of rank $\kappa$ on a vector space
of dimension $n$.
Assume that $n\geq2^\kappa$ and that  $\kappa=4,5$. Then:
\begin{enumerate}
\item Suppose that $\mathcal{V}$ is complex Osserman. Then 
there exists a reparametrization $\tilde{\mathcal{F}}=\{\tilde J_1,...,\tilde J_\kappa\}$ of
the family $\mathcal{F}=\{J_1,...,J_\kappa\}$
so that $J=\tilde J_1\tilde J_2$ and so that 
$R=\tilde c_1R_{\tilde J_1}+\dots+\tilde c_\kappa R_{\tilde J_\kappa}$.
\item If $\kappa=5$, then $\mathcal{V}$ is not complex Osserman.
\item If $\kappa=4$, then $\mathcal{V}$ is not complex Osserman.
\end{enumerate}
\end{lemma}

\begin{proof} Since $\kappa=4$ or $\kappa=5$ we have
$2\kappa+5<16\leq n$. Thus Lemma \ref{lem-4.1} implies
$Jx \in\operatorname{Span}_{i\ne j}\{J_iJ_jx\}$ for all $x\in
S(V,\langle\cdot,\cdot\rangle)$.
One can show there exists $x,y\in V$ so $\{J_jJ_kx,J_jJ_ky\}_{j<k}$ is an orthonormal set of
$\kappa(\kappa-1)$ linearly independent vectors. Thus the argument used to establish Lemma \ref{lem-3.2}
proves that
\[
J=\sum_{i=1}^{\kappa-1}\sum_{j=i+1}^\kappa a_{ij} J_iJ_j\,.
\]
One can now show that there exists a suitable reparametrization; as the argument is
straightforward, if a bit lengthy, we shall omit the details.

Suppose that $\kappa=5$. By Assertion (1), we may suppose
that $J=J_1J_2$. As noted above, there
exists $x\in S(V,\langle\cdot,\cdot\rangle)$ such that
$\{J_iJ_jx\}_{i<j}$ is an orthonormal set and, thus,
$\{J_1x,J_2x,J_3x,J_4x,J_5x,J_1J_2J_3x,J_1J_2J_4x,J_1J_2J_5x\}$ is
also an orthonormal set. The\-refore
\[
\jp y=\left\{\begin{array}{ll}
                3(c_1+c_2)y&\mbox{if}\;y\in \operatorname{Span}\{J_1x,J_2x\},\\
                3c_3y&\mbox{if}\;y\in \operatorname{Span}\{J_3x,J_1J_2J_3x\},\\
                3c_4y&\mbox{if}\;y\in \operatorname{Span}\{J_4x,J_1J_2J_4x\},\\
                3c_5y&\mbox{if}\;y\in \operatorname{Span}\{J_5x,J_1J_2J_5x\},\\
                0&\mbox{otherwise}.
                \end{array}\right.
\]
Note that $\operatorname{Rank}\{\jp y\}>4$. Hence, by Lemma \ref{thm-2.4}, $R$ is
not complex Osserman. Assertion (2) now follows.

Finally suppose $\kappa=4$. Again, we may suppose
$J=J_1J_2$. Since $(J_1J_2J_3)^2=\operatorname{id}$, there exists $x\in
S(V,\langle\cdot,\cdot\rangle)$ such that $J_1J_2J_3x=\pm x$, and
hence 
$$\{x,J_1x,J_2x,J_3x,J_4x,J_1J_2J_4x\}$$
is an orthonormal
set. Note that
\[
\jp y=\left\{\begin{array}{ll}
                3c_3y&\mbox{if}\;y\in \operatorname{Span}\{x,J_3x\},\\
                3(c_1+c_2)y&\mbox{if}\;y\in \operatorname{Span}\{J_1x,J_2x\},\\
                3c_4y&\mbox{if}\;y\in \operatorname{Span}\{J_4x,J_1J_2J_4x\},\\
                0&\mbox{if}\;y\bot\,\operatorname{Span}
                \{x,J_3x,J_1x,J_2x,J_4x,J_1J_2J_4x\}.
\end{array}\right.
\]
Now, since $(J_1J_2J_3J_4)^2=\operatorname{id}$, there exists $y\in
S(V,\langle\cdot,\cdot\rangle)$ such that $J_1J_2J_3J_4y=\pm y$
and
\[
\jp y=\left\{\begin{array}{ll}
                3(c_1+c_2)y&\mbox{if}\;y\in \operatorname{Span}\{J_1x,J_2x\},\\
                3(c_3+c_4)y&\mbox{if}\;y\in \operatorname{Span}\{J_3x,J_4x\},\\
                0&\mbox{if}\;y\bot\, \operatorname{Span}\{J_1x,J_2x,J_3x,J_4x\}.
                \end{array}\right.
\]
Since the eigenvalues are different, $R$ is not complex Osserman.
\end{proof}

\subsection{Curvature given by a Clifford family with $c_0\neq 0$.}\label{sect-4.2}
This section is devoted to the proof of Assertion
(2) of Theorem \ref{thm-1.3}. Although there is some
parallelism between cases $c_0=0$ and $c_0\neq 0$, the approach we
follow now is slightly different. However, in the interests of brevity, we
will refer to arguments in Section \ref{sect-4.1} whenever possible. We begin
by studying a reduced complex Jacobi operator where the effect of $c_0$ has been
normalized.

\begin{lemma}\label{lem-4.4}
Let $\mathcal{V}:=(V,\g\cdot\cdot,J,R=c_0R_0+c_1R_{J_1}+\dots+c_\kappa
R_{J_\kappa})$ be a complex Osserman model where
$\{J_1,...,J_\kappa\}$ is a Clifford family of rank $\kappa$ on a vector space
of dimension $n$. Assume that $\kappa\geq 4$.
If $4\leq \kappa \leq 7$, assume that $n\geq 2^\kappa$. If $\kappa\geq 8$,
assume that $n\geq\kappa(\kappa-1)$. Let
$\tilde{\mathcal{J}}_R(\pi_x)=\mathcal{J}_R(\pi_x)-2c_0\operatorname{id}$. Then:
\begin{enumerate}
\item $\operatorname{Rank}\{\tilde{\mathcal{J}}_R(\pi_x)\}\leq 4$. \item $Jx\in
\operatorname{Span}\{J_ix,J_jJ_kx\}_{i,j<k}$ for all $x\in V$. \item If
$\kappa\geq 6$, then $Jx\in \operatorname{Span}\{J_iJ_jx\}_{i\leq 6}$ for all
$x\in V$.
\item $\kappa\le5$.
\end{enumerate}
\end{lemma}
\begin{proof}
We use Equation (\ref{eqn-1.g}) to see that:
\[
\tilde{\mathcal{J}}_R(\pi_x)y=-c_0\g{y}{x}x-c_0\g{y}{Jx}Jx+3\sum
c_i(\g{y}{J_ix}J_ix+\g{y}{J_iJx}J_iJx)\,.
\]
Consequently $\operatorname{Rank}\{\tilde{\mathcal{J}}_R(\pi_x)\}\leq 2\kappa +2$ and $0$ is
an eigenvalue with multiplicity at least $n-2\kappa-2$. Since $n-2\kappa-2>4$ and as we
have simply shifted the spectrum,
Theorem \ref{thm-2.4} may be used to derive Assertion (1).

To prove Assertion (2), we compute that:
\begin{eqnarray}\label{eqn-4.b}
&&\mathcal{\tilde J}_R(\pi_x)x=-c_0x+\textstyle\sum_i3c_i\langle x,J_iJx\rangle J_iJx,\nonumber\\
&&\mathcal{\tilde J}_R(\pi_x)Jx=-c_0Jx+\textstyle\sum_i3c_i\langle Jx,J_ix\rangle J_ix,\\
&&\mathcal{\tilde J}_R(\pi_x)J_ix=-c_0\langle J_ix,Jx\rangle
Jx+3c_iJ_ix+\textstyle\sum_j3c_j\langle J_ix,J_jJx\rangle J_jJx\,.
\nonumber
\end{eqnarray}
Define:
\begin{eqnarray*}
&&M:=\operatorname{diag}(-c_0,3c_1,...,3c_\kappa),\\
&&U(x):=\operatorname{Span}\{x,J_1x,...,J_\kappa x\},\\
&&V(x):=\operatorname{Span}\{Jx,J_1Jx,...,J_\kappa Jx\},\\
&&W(x):=U(x)+V(x)\,.
\end{eqnarray*}
Let $\rho$ denote projection on $W(x)/V(x)$. We then have that
$\rho\tilde{\mathcal{J}}(\pi_x)=\rho M$ on $U(x)$. As $M$ is
invertible, the following inequalities hold:
\begin{eqnarray*}
&&\dim\{\rho U(x)\}=\dim\{\rho\mathcal{\tilde J}_R(\pi_x)U(x)\}\le 4,\\
&&\dim\{W(x)\}\le4+\kappa+1,\\
&&\dim\{U(x)\cap V(x)\}\ge\kappa+1+\kappa+1-\kappa-5=\kappa-3>0\,.
\end{eqnarray*}
Therefore, there exists a non-trivial relationship
\[
(a_0+a_1J_1+...+a_\kappa J_\kappa)Jx=(b_0+b_1J_1+...+b_\kappa
J_\kappa)x\,.
\]
We invert this relationship by multiplying by $(a_0-a_1J_1-...-a_\kappa J_\kappa)$. Since
$Jx\perp x$, we may conclude that $Jx\in \operatorname{Span}\{J_ix,J_jJ_kx\}$ and
establish Assertion (2).

If $\kappa\geq 6$, then we can derive a stronger result. We estimate that:
\begin{eqnarray*}
&&\dim\big\{\{\operatorname{Span}\{J_1x,...,J_6x\}\cap\operatorname{Span}\{J_1Jx,...,J_\kappa
Jx\}\big\}\\
&&\ge 6+\kappa-\dim(W)\ge6+\kappa-\kappa-5>0\,.
\end{eqnarray*}
Assertion (3) now follows using a similar argument to that used to establish Assertion (2).

To establish Assertion (4), we assume to the contrary that $\kappa\ge6$ and argue for a
contradiction. By Assertion (3), we have that 
$Jx\in \operatorname{Span}\{J_iJ_jx\}_{i\leq 6,j\ne i}$. The argument used to establish 
Lemma \ref{lem-4.2} shows that $\kappa\leq 7$. Thus we have that $\kappa=6$ or $\kappa=7$. 
Since $n\geq 2\kappa(\kappa-1)$, Corollary \ref{cor-3.3} and
Assertion (3) show that $J\in \operatorname{Span}\{J_iJ_j\}$. One may show
there exists $x\in V$ such that $x\bot J_iJ_jJ_kx$ for
any $i,j,k$ and such that $J_1J_2x\bot
\operatorname{Span}\{J_iJ_jx\}_{(i,j)\neq (1,2)}$. Thus, since $Jx\bot J_ix$ for
this specific $x$, Equation (\ref{eqn-4.b})
yields:
\begin{eqnarray*}
&&\mathcal{\tilde J}_R(\pi_x)x=-c_0x,\qquad\mathcal{\tilde J}_R(\pi_x)Jx=-c_0Jx,\\
&&\mathcal{\tilde J}_R(\pi_x)J_ix=3c_iJ_ix+\displaystyle\sum_{j=1}^\kappa3c_j\langle
J_ix,J_jJx\rangle J_jJx\,.
\end{eqnarray*}
Hence the subspace $\operatorname{Span}\{x,Jx\}$ is invariant under
$\mathcal{\tilde J}(\pi_x)$. We clear the previous notation. By applying the argument used to prove
Assertion (2) to the sets
\begin{eqnarray*}
&&U(x):=\operatorname{Span}\{J_1x,...,J_\kappa x\},\\
&&V(x):=\operatorname{Span}\{J_1Jx,..., J_\kappa Jx\},\\
&&W(x):=U(x)+V(x),
\end{eqnarray*}
we obtain $Jx=\sum_{i\leq 3,i<j} a_{ij} J_iJ_jx$. Thus in particular $a_{45}=0$.
Since the coefficients $a_{ij}$
were universal and independent of $x$, we can permute the indices to see that 
$a_{ij}=0$ for all $i<j$, which is impossible.
\end{proof}

It remains to show that a Clifford family of rank $\kappa=4$ or $\kappa=5$ can
not give a complex Osserman model. As in the
case $c_0=0$ these ranks are treated independently. However, the present
situation is a bit more difficult. We present sketch of
proofs describing the main ideas involved; full details are available from the authors
upon request but are omitted here in the interests of brevity.

\begin{lemma}\label{lem-4.5}
Let $\mathcal{V}:=(V,\g\cdot\cdot,J,R=c_0R_0+c_1R_{J_1}+\dots+c_\kappa
R_{J_\kappa})$ be a complex Osserman model where
$\{J_1,...,J_\kappa\}$ is a Clifford family of rank $\kappa=4$ or $\kappa=5$ on a vector space
of dimension $n\ge32$. 
\begin{enumerate}
\item If $\kappa=5$, then $\mathcal{V}$ is not complex Osserman.
\item If $\kappa=4$, then $\mathcal{V}$ is not complex Osserman.
\end{enumerate}
\end{lemma}

\begin{proof} Suppose that $\kappa=5$, that $n\ge32$, and that $\mathcal{V}$ is complex Osserman.
We argue for a contradiction. Using similar techniques to those which were used to prove
Lemma \ref{lem-4.4}, one shows that $J\notin
\operatorname{Span}\{J_iJ_j\}_{i\ne j}$. Consider the set
\[
C:=\{x\in V: Jx\in \operatorname{Span}\{J_ix\}\}\,.
\]
One shows that $C$ is a closed
nowhere dense set. So, working in the complementary set $C^c$ and
using similar arguments to those which were used to prove
Lemma \ref{lem-4.1} applied to the sets
\begin{eqnarray*}
&&U(x):=\operatorname{Span}\{J_1x,...,J_5x,Jx\},\\
&&V(x):=\operatorname{Span}\{J_1Jx,...,J_5Jx\},\\
&&W(x):=U(x)+V(x),
\end{eqnarray*}
one shows that $Jx\in \operatorname{Span}\{J_iJ_jx\}_{i\ne j}$ and, therefore, $J\in
\operatorname{Span}\{J_iJ_j\}_{i\ne j}$, which is false. This proves Assertion (1).

Suppose that $\kappa=4$. By Lemma \ref{lem-4.4} we know that $Jx\in
\operatorname{Span}\{J_ix,J_jJ_kx\}_{j<k}$ for all $x\in V$. Since $n\ge32$,
one can show that there exists $x,y\in S(V,\g\cdot\cdot)$ so that
$\{J_ix,J_{jk}x,J_iy,J_{jk}y\}_{j<k}$ is an orthonormal set. The argument given to establish
Lemma \ref{lem-3.2} (3) then shows there exist constants $a_i$ and
$a_{jk}$ so that
\[
J=\sum_{i=1}^4a_iJ_i+\sum_{j<k}a_{jk}J_jJ_k\,.
\]
The compatibility between $J$ and $R$ shows that the constants $a_i$ vanish so
\[
J=\sum_{i<j}a_{ij}J_iJ_j\,.
\]
In this situation one may reparametrize the Clifford family so
$J=\tilde J_1\tilde J_2$. A straightforward calculation now shows
$\operatorname{Rank}\{\tilde{\mathcal{J}}_{\pi_x}\}\geq 6$, which contradicts
Theorem \ref{thm-2.4}.\end{proof}

\section{Classification for Clifford families of lower rank}\label{sect-5}

In this section we prove Theorem \ref{thm-1.4} by studying complex Osserman models which are given by
Clifford families of rank
$\kappa$ for
$0\le\kappa\le 3$. Section \ref{sect-5.1} deals with the case $\kappa=0$, Section \ref{sect-5.2}
deals with $\kappa=1$, and Section \ref{sect-5.3} deals with $\kappa=2$. We
shall omit much of the analysis when discussing the case $\kappa=3$ in Section \ref{sect-5.4} in the
interests of brevity as it is  similar to the other cases; again, details are available upon request from
the authors.
 Throughout Section \ref{sect-5}, we suppose that $R=c_0R_0+c_1R_{J_1}+...+c_\kappa R_{J_\kappa}$.

\subsection{Clifford families of rank $0$}\label{sect-5.1}
Let $\mathcal{V}=(V,\g{\cdot}{\cdot},J,R=c_0R_0)$. Then we have
\[
\ja(\pi_x)y=c_0(2y-\g{y}{x}x-\g{y}{Jx}Jx).
\]
Hence $J\ja(\pi_x)=\ja(\pi_x)J$ and the eigenvalues are
$(c_0,2c_0)$ with multiplicities $(2,n-2)$ for any $x\in S(V,\langle\cdot,\cdot\rangle)$. Consequently,
$R$ is complex Osserman.

\subsection{Clifford families of rank $1$}\label{sect-5.2}

We have that:

\begin{lemma}\label{RJ1}
Let $J$ and $J_1$ be Hermitian almost complex structures on $(V,\langle\cdot,\cdot\rangle)$.
\begin{enumerate}
\item 
Let $\mathcal{V}:=(V,\g\cdot\cdot,J,R=c_1R_{J_1})$ where $c_1\ne0$.
The following assertions are equivalent:
\begin{enumerate}
\item $R$ and $J$ are compatible.
\item $JJ_1=J_1J$ or $JJ_1=-J_1J$.
\item $\mathcal{V}$ is complex Osserman.
\end{enumerate}
\item Let $\mathcal{V}=(V,\g{\cdot}{\cdot},J,R=c_0R_0+c_1R_{J_1})$ where $c_0c_1\ne0$. Then
$\mathcal{V}$ is complex Osserman if and only if $J=\pm J_1$ or
$JJ_1=-J_1J$.
\end{enumerate}
\end{lemma}

\begin{proof} Suppose that $\mathcal{V}:=(V,\g\cdot\cdot,J,R=c_1R_{J_1})$
and that $J$ and $R$ are compatible. By 
Equation (\ref{eqn-1.g}),
\[
\ja(\pi_x)y=3\g{y}{J_1x}J_1x+3\g{y}{J_1Jx}{J_1Jx}\,.
\]
Hence $\operatorname{Range}\{\ja(\pi_x)\}=\operatorname{Span}\{J_1x,J_1Jx\}$ and,
since $J$ and $R$ are compatible, we have
$J(\operatorname{Span}\{J_1x,J_1Jx\})\subset
\operatorname{Span}\{J_1x,J_1Jx\}$. Since $JJ_1x\bot J_1x$,
necessarily $JJ_1x=\epsilon_x J_1Jx$, where $\epsilon_x=\pm 1$. By
continuity, since $S(V,\langle\cdot,\cdot\rangle)$ is connected,
$\epsilon_x$ is constant. Then $JJ_1=J_1J$ or $JJ_1=-J_1J$. If this condition holds,
then it is easily verified that $\mathcal{V}$ is complex Osserman. Finally, if $\mathcal{V}$
is complex Osserman, then necessarily $R$ and $J$ are compatible. Assertion (1) now follows.

Next suppose that $\mathcal{V}:=(V,\g\cdot\cdot,J,R=c_0R_0+c_1R_{J_1})$
is a complex Osserman model where $c_0\ne0$ and $c_1\ne0$. 
Since $R$ and $J$ are compatible and since $R_0$
and $J$ are compatible, $R_{J_1}$ and $J$ are compatible as well.
Thus by Assertion (1),
$JJ_1=J_1J$ or $JJ_1=-J_1J$. We now
show that $JJ_1=J_1J$ implies $J=\pm J_1$. 
We suppose to the contrary that $J\ne\pm J_1$
and argue for a contradiction. Because $(JJ_1)^2=\operatorname{id}$, we can use
$JJ_1$ to define a $\mathbb{Z}_2$ grading on $V$ by decomposing $V=V_+\oplus V_-$ where $J=\pm J_1$
on $V_\pm$.

Let $x_\pm\in S(V_\pm)$ and let $x_0=(x_++x_-)/\sqrt{2}$. Then one has that:
\[
\begin{array}{rcl}
\ja(\pi_{x_+})y&=&\left\{\begin{array}{rl}
                            (c_0+3c_1)y&\phantom{..}\mbox{if}\; y\in \operatorname{Span}\{x_+,Jx_+\},\\
                            2c_0y &\phantom{..}\mbox{if}\; y\bot\,
                            \operatorname{Span}\{x_+,Jx_+\},
\end{array}\right.\\
\ja(\pi_{x_0})y&=&\left\{\begin{array}{rl}
                            c_0y &\mbox{if}\; y\in \operatorname{Span}\{x_+,Jx_0\},\\
                            (2c_0+3c_1)y&\mbox{if}\; y\in \operatorname{Span}\{J_1x_0,J_1Jx_0\},\\
                            2c_0y &\mbox{if}\; y\bot\,
                            \operatorname{Span}\{x_0,Jx_0,J_1x_0,JJ_1x_0\}.
\end{array}\right.
\end{array}
\]
This shows that the eigenvalues of $\ja(\pi_{x_+})$ are $(c_0+3c_1,2c_0)$
with multiplicities $(2,n-2)$ (if $3c_1=c_0$ then $2c_0$ has
multiplicity $n$). Furthermore, the eigenvalues of $\ja(\pi_{x_0})$ are
$(c_0,2c_0+3c_1,2c_0)$ with multiplicities $(2,2,n-4)$. So the
eigenvalues are different in both cases. This contradiction shows that if $JJ_1=J_1J$,
then $J=\pm J_1$.

Conversely, if $JJ_1=-J_1J$ or if $J=\pm J_1$, then a straightforward calculation
shows $\mathcal{V}$ is complex Osserman.
\end{proof}

\subsection{Clifford families of rank $2$}\label{sect-5.3}
We first suppose that $c_0=0$.

\begin{lemma}\label{lem-5.2}
Let $J$ be an Hermitian almost complex structure and let $\{J_1,J_2\}$ be a
Clifford family on $(V,\langle\cdot,\cdot\rangle)$.
Let $\mathcal{V}:=(V,\g{{\cdot}}{{\cdot}},J,R=c_1R_{J_1}+c_2R_{J_2})$ be complex Osserman.
If $x$ is a unit vector, set $\alpha(x):=\g{J_1J_2x}{Jx}$. Then:
\begin{enumerate}
\item $\alpha(x)$ is constant on
$S(V,\langle\cdot,\cdot\rangle)$.
\item Either $\alpha=0$, or $\alpha=1$, or $\alpha=-1$.
\item Suppose that $\alpha=\pm 1$. Then $J=\pm J_1J_2$ and $\operatorname{Rank}\{\mathcal{J}_R(\pi_x)\}=2$.
\item Suppose that $\alpha=0$. Then $\operatorname{Rank}\{\mathcal{J}_R(\pi_x)\}=4$. Furthermore:
\begin{enumerate}
\item if $c_1\neq c_2$ then $J J_1=J_1 J$ and $J J_2=-J_2 J$ or $J
J_1=-J_1 J$ and $J J_2=J_2 J$.
\item if $c_1=c_2$ then there
exists a reparametrization $\{\tilde J_1,\tilde J_2\}$ of
$\{J_1,J_2\}$ so that $R=c_1 R_{\tilde J_1}+c_2 R_{\tilde J_2}$,
$J \tilde J_1=\tilde J_1 J$ and $J \tilde J_2=-\tilde J_2 J$.
\end{enumerate}
\end{enumerate}
\end{lemma}

\begin{proof}
Since $\mathcal{V}$ is complex Osserman, Equation (\ref{eqn-1.g})
shows that 
$$
\operatorname{Range}\{\ja(\pi_x)\}\subset
\operatorname{Span}\{J_1x,J_1Jx,J_2x,J_2Jx\}\,.
$$
Consequently,
\begin{eqnarray*}
&&\jp J_1x=3c_1J_1x+3\alpha(x) c_2J_2Jx,\\
&&\jp J_2Jx=3\alpha(x)c_1J_1x+3c_2J_2Jx,\\
&&\jp J_1Jx=3c_1J_1Jx-3\alpha(x)c_2J_2x,\\
&&\jp J_2x=-3\alpha(x)c_1J_1Jx+3c_2J_2x\,.
\end{eqnarray*}
Thus $V_1(x):=\operatorname{Span}\{J_1x,J_2Jx\}$ and
$V_2(x):=\operatorname{Span}\{J_2x,J_1Jx\}$ are $\jp$ invariant subspaces. 
Note that $J(V_1(x))=V_2(x)$, that
$V_1(x)\bot V_2(x)$, and that 
$$\operatorname{Range}(\jp)=V_1(x)\oplus V_2(x)\,.$$

If $\alpha(\bar x)=\pm 1$ for some $\bar x\in
S(V,\langle\cdot,\cdot\rangle)$, then $\operatorname{Rank}\{\mathcal{J}_R(\pi_{\bar x})\}=2$. Since
$R$ is complex Osserman, $\jp$ has constant rank. In such a
case we get $\alpha(x)=\pm 1$ for all $x\in
S(V,\langle\cdot,\cdot\rangle)$. On the other hand if $\alpha(x)\neq \pm 1$,
then
\begin{eqnarray*}
&&\jp|_{V_1(x)}=\left(\begin{array}{cc}
                                        3c_1&3\alpha(x)c_1\\
                                        3\alpha(x)c_2&3c_2
\end{array}\right)\!\!,\quad\text{and}\\
&&\jp|_{V_2(x)}=\left(\begin{array}{cc}
                                        3c_1&-3\alpha(x)c_1\\
                                        -3\alpha(x)c_2&3c_2
\end{array}\right)\!.
\end{eqnarray*}
Consequently, $\det\{\jp|_{V_1(x)+V_2(x)}\}=(9c_1c_2(1-\alpha(x)^2))^2$. Since the eigenvalues
of $\ja(\cdot)$ are constant, the determinant of $\ja(\cdot)$ is constant and consequently
$\alpha(x)$ does not depend on $x$. This establishes Assertion (1). The proof of 
Assertion (2) is a bit technical and is omitted in the interests of brevity.
It relies on the fact that $J$ preserves the eigenspaces of $\jp$; details
are available from the authors.

The possible values of $\operatorname{Rank}\{\jp\}$ are $2$ and $4$, which
correspond to $\alpha=\pm 1$ or $\alpha\neq \pm 1$, respectively.
If $\alpha=\pm 1$, then $J=\pm J_1J_2$ since $Jx$
and $J_1J_2x$ are unit vectors. Assertion (3) now follows.

On the other hand, if $\alpha=0$
then, by polarizing the identity $\g{J_1J_2x}{Jx}=0$, we see that
$\g{J_1J_2x}{Jy}+\g{J_1J_2y}{Jx}=0$ and consequently $J_1J_2J+JJ_1J_2=0$.
Furthermore, $\{J_1x,J_1Jx,J_2x,J_2Jx\}$ is an
orthonormal set for any $x\in S(V,\g\cdot\cdot)$.

Suppose that $c_1\neq c_2$ and that $\ja$ has three different
eigenvalues $(0,3c_1,3c_2)$. As $J$ preserves the
eigenspaces of $\jp$, $J$ preserves the spaces $\operatorname{Span}\{J_1x,J_1 Jx\}$ and
$\operatorname{Span}\{J_2x,J_2Jx\}$. Consequently, $J J_1=\pm J_1 J$ and $J J_2=\pm J_2 J$.
Since one has that $JJ_1J_2+J_1J_2J=0$,
the only possibilities are $J J_1=J_1 J$ and $J J_2=-J_2 J$ or
$J J_1=-J_1 J$ and $J J_2=J_2 J$.

Suppose that $c_1=c_2$. In such a case there are only two distinct
eigenvalues for $\jp$ and
$\operatorname{Range}\{\jp\}=\operatorname{Span}\{J_1x,J_2x,J_1Jx,J_2Jx\}$ is a
$4$-dimensional eigenspace. Since $J$ preserves this eigenspace
and $J_1Jx\bot J_1x,J_2x$ we have
\[
JJ_1x=\g{JJ_1x}{J_1Jx}J_1Jx+\g{JJ_1x}{J_2Jx}J_2Jx\,.
\]
 Set $\Theta_1=JJ_1$ and $\Theta_2=JJ_2$, then $\g{\Theta_1^2x}{x}^2+\g{\Theta_2\Theta_1x}{x}^2=1$.
 Also note that
\begin{eqnarray*}
&&\Theta_1 \Theta_1^*=JJ_1J_1J=\operatorname{id},\Theta_2 \Theta_2^*=JJ_2J_2J=\operatorname{id},\\
&&\Theta_1\Theta_2^*+\Theta_2\Theta_1^*=JJ_1J_2J+JJ_2J_1J=0,\\
&&\Theta_1\Theta_2=JJ_1JJ_2=JJ_1JJ_1J_2J_1=-JJ_1J_1J_2JJ_1=\Theta_2\Theta_1\,.
\end{eqnarray*}
Consequently, $\Theta_1$ and $\Theta_2$ are commuting orthogonal
maps. Let
\[
V=V_+\oplus V_-\oplus V_1\oplus\dots\oplus V_k
\]
be a skew-diagonalization of $\Theta_1$, such that $\Theta_1=\pm
\operatorname{id}$ on $V_\pm$ and $\Theta_1$ is a rotation through an angle
$\theta_i$, $0<\theta_i<\pi$, on $V_i$. After some technical fuss,
one may show that there is a reparametrization $\{\tilde
J_1,\tilde J_2\}$ such that the previous decomposition is reduced
to $V=V_+\oplus V_-$ and hence $J \tilde J_1=\tilde J_1 J$. Also,
since $JJ_1J_2=-J_1J_2J$ as noted above,
$J\tilde J_2=-\tilde J_2J$.
\end{proof}

We complete the proof of Theorem \ref{thm-1.4} (3) by
studying models with $c_0\neq 0$. First we
establish the following consequence of the compatibility between
$J$ and $R$ for a Clifford family of rank at most $3$.

\begin{lemma}\label{lem-5.3}
Let $R=c_1R_{J_1}+c_2R_{J_2}+c_3R_{J_3}$ be an algebraic curvature
tensor given by a Clifford family of rank $3$. Suppose $R$ is compatible with
an Hermitian almost complex structure $J$. If
$Jx=(a_1J_1+a_2J_2+a_3J_3)x$ for all $x\in V$, then
$(c_i-c_j)a_ia_j=0$ for $i\ne j$.
\end{lemma}
\begin{proof}
Compute
\begin{eqnarray*}
&&J R(x,Jx)x=c_0x-3c_1a_1JJ_1x-3c_2a_2JJ_2x-3c_3a_3JJ_3x,\\
&&R(x,Jx)Jx=c_0x-3c_1a_1J_1Jx-3c_2a_2J_2Jx-3c_3a_3J_3Jx.
\end{eqnarray*}
Now, since $R$ and $J$ are compatible, $J R(x,Jx)x=R(x,Jx)Jx$ so
\[
(c_1-c_2)a_1a_2J_1J_2x+(c_1-c_3)a_1a_3J_1J_3x+(c_2-c_3)a_2a_3J_2J_3x=0.
\]
Since $\{J_1J_2x,J_1J_3x,J_2J_3x\}$  is an orthogonal set, the
desired equalities follow.
\end{proof}

\begin{lemma}\label{lem-5.4}
Let
$\mathcal{V}=(V,\g{{\cdot}}{{\cdot}},J,R=c_0R_0+c_1R_{J_1}+c_2R_{J_2})$
be complex Osserman. If $\dim\{V\}\geq 12$, then there exists a
reparametrization $\{\tilde J_1,\tilde J_2\}$ of $\{J_1,J_2\}$
such that $R=c_0R_0+\tilde c_1R_{\tilde J_1}+\tilde c_2R_{\tilde J_2}$ and
either $J=\tilde J_1$ or $J=\tilde J_1\tilde J_2$.
\end{lemma}
\begin{proof}
Let $\tjp=\jp-2c_0 \operatorname{id}$ be the reduced complex Jacobi operator. As $\jp$ is complex
Osserman,
$\tjp$ has rank at most $4$. Let $\alpha(x):=\g{J_1J_2x}{Jx}$,
$\alpha_1(x):=\g{J_1x}{Jx}$ and $\alpha_2(x):=\g{J_2x}{Jx}$. Then
\begin{eqnarray*}
\tjp x&=& -c_0 x-3c_1\alpha_1(x) J_1Jx-3c_2\alpha_2(x)J_2Jx,\\
\tjp Jx&=& -c_0 Jx+3c_1\alpha_1(x) J_1x+3c_2\alpha_2(x)J_2x,\\
\tjp J_1x&=& -c_0\alpha_1(x)Jx+3c_1 J_1x+3c_2\alpha(x) J_2Jx,\\
\tjp J_2x&=& -c_0\alpha_2(x)Jx-3c_1\alpha(x) J_1Jx+3c_2J_2x,\\
\tjp J_1Jx&=& c_0\alpha_1(x) x+3c_1 J_1Jx-3c_2\alpha(x) J_2x,\\
\tjp J_2Jx&=& c_0\alpha_2(x) x+3c_1\alpha(x) J_1x+3c_2 J_2Jx.
\end{eqnarray*}
Consider the subspace $W(x):=\operatorname{Span}\{x,J_1x,J_2x,Jx,J_1Jx,J_2Jx\}$
and notice that $\operatorname{Range}\{\tjp\}\subset W(x)$. We wish to show
that $\dim\;W(x)<6$. On the contrary, suppose $\dim\{W(x)\}=6$. From
the previous calculations we get the matrix associated to
$\tjp_{|W(x)}$ and compute:
$$\det(\tjp_{|W(x)})=3^4c_0^2c_1^2c_2^2(-1+\alpha(x)^2+\alpha_1(x)^2+\alpha_2(x)^2)^2\,.$$
Since $\dim\{V\}\geq 12$ we apply Theorem \ref{thm-2.4}
to get $\det(\tjp_{|W(x)})=0$ and hence
$\alpha^2+\alpha_1^2+\alpha_2^2=1$. Since $\alpha(x)$,
$\alpha_1(x)$ and $\alpha_2(x)$ are the Fourier coefficients of
$Jx$ with respect to $\{J_1J_2x,J_1x,J_2x\}$, we get $Jx=\alpha(x)
J_1J_2x+\alpha_1(x) J_1x+\alpha_2(x)J_2x$ which contradicts
the assumption that $\dim\{W(x)\}=6$.

Hence $\dim\{W(x)\}\leq 5$ and
$\operatorname{Span}\{x,J_1x,J_2x\}\cap
\operatorname{Span}\{Jx,J_1Jx,J_2Jx\}$ is non-trivial. Moreover,
there exists a unit vector $(\rho_0,\rho_1,\rho_2)\in \mathbb{R}^3$
such that 
$$(\rho_0+\rho_1J_1+\rho_2J_2)Jx\in
\operatorname{Span}\{x,J_1x,J_2x\}\,.$$
Let $\{J_1,J_2,J_1J_2\}$ give $V$ a quaternion structure $\mathbb{H}$. As
$Jx\in \mathbb{H}x$,
\[
Jx=a_1(x)J_1x+a_2(x)J_2x+a_3(x)J_3x.
\]
The following argument shows that $a_i(\cdot)$ are constant
functions in $S(V,\langle\cdot,\cdot\rangle)$. Let $x,y\in S(V,\langle\cdot,\cdot\rangle)$. Since
$\dim\{\mathbb{H}x+\mathbb{H}y\}\leq 8$, there exists $z\in S(V,\langle\cdot,\cdot\rangle)$
such that $z\bot \mathbb{H}x,\mathbb{H}y$. Then $\mathbb{H}x\bot
\mathbb{H}z$ and for $w:=\frac{1}{\sqrt{2}}(x+z)$ we have:
\[
J(w)=\frac{1}{\sqrt{2}}\sum_i
a_i(w)J_i(x+z)=\frac{1}{\sqrt{2}}\sum_i(a_i(x)J_i(x)+a_i(z)J_i(z)).
\]
which implies that $a_i(x)=a_i(w)=a_i(z)$. Similarly, $a_i(y)=a_i(z)$.

Therefore $J=a_1J_1+a_2J_2+a_3J_1J_2$. By Lemma
\ref{lem-5.3} with $c_3=0$, we have:
$$(c_1-c_2)a_1a_2=c_2a_2a_3=c_1a_1a_3=0\,.$$
Then either $J=\pm J_3$ or $J=a_1J_1+a_2J_2$ and we may
reparametrize $\{J_1,J_2\}$ by $\{\tilde J_1,\tilde J_2\}$ so that
$J=\tilde J_1$.
\end{proof}

\subsection{Clifford families of rank $3$}\label{sect-5.4}

Let $\{J_1,J_2,J_3\}$ be a Clifford family on $V$. The
dual structure, which is always a quaternion structure, is given by
\[
\{J_1^*:=J_2J_3,J_2^*:=J_3J_1,J_3^*:=J_1J_2\}\,.
\]

We use this structure to establish Assertion (4) of  Theorem \ref{thm-1.4}; in the interest of brevity
we shall simply outline the proof rather than giving full details. Let
$\mathcal{V}=(V,\g\cdot\cdot,J,R=c_0R_0+c_1R_{J_1}+c_2R_{J_2}+c_3R_{J_3})$
be complex Osserman, where $c_0$ may be $0$. Then one has the
following:
\begin{enumerate}
\item If $J=a_1J_1+a_2J_2+a_3J_3$, then there exists a
reparametrization $\{\tilde J_1,\tilde J_2,\tilde J_3\}$ so that
$R=c_0R_0+\tilde c_1R_{\tilde J_1}+\tilde c_2R_{\tilde J_2}+\tilde
c_3R_{\tilde J_3}$ and $J=\tilde J_1$. \item Suppose $J_1J_2\neq
J_3$.
\begin{enumerate}
\item Then $J\neq J_1$. Furthermore, if $c_0\neq 0$, then $J\neq
J_2J_3$. \item Suppose that $Jx\in
\operatorname{Span}\{J_1x,J_2x,J_3x,J_1^*x,J_2^*x,J_3^*x\}$ for some element $x\in V$ with
$x=(x_++x_-)/\sqrt{2}$ where
$J_1J_2x_\pm=J_3x_\pm$. Then $c_0=0$, and there is a
reparametrization $\{\tilde J_1,\tilde J_2,\tilde J_3\}$ such that one has
 $R=\tilde c_1R_{\tilde J_1}+\tilde
c_2R_{\tilde J_2}+\tilde c_3R_{\tilde J_3}$ and $J=\tilde
J_2\tilde J_3$.
\end{enumerate}
\end{enumerate}

The classification in Theorem \ref{thm-1.4} (4) follows from these
observations and from a careful analysis of the rank of
the matrix associated to $\jp$. The technique is similar to that
developed in Lemma \ref{lem-5.4}.

\section*{Acknowledgments}
The research of M. Brozos-V\'azquez was partially supported by
projects BFM 2003-02949 and PGIDIT04PXIC20701PN (Spain) and by the Max Planck
Institute for Mathematics in the Sciences (Leipzig, Germany). The
research of P. Gilkey was partially supported by the Max Planck
Institute for Mathematics in the Sciences (Leipzig, Germany),
by DFG-project 158/4-4, and by the Spectral Theory and Partial Differential Equations Program of the
Newton Institute (Cambridge, U.K.).

It is a pleasure for both authors to acknowledge helpful
conversations with Prof. E. Garc{\'\i}a-R{\'\i}o.


\begin{thebibliography}{}

\bibitem{Adams} Adams, J., `Vector fields on spheres', {\it Ann. of Math.} {\bf
75} (1962), 603--632.

\bibitem{ABS}
Atiyah, M. F., Bott, R., and Shapiro, A.,
`Clifford Modules', {\it Topology} {\bf 3} suppl. 1 (1964), 3--38.

\bibitem{Chi} Chi, Q.-S.,
`A curvature characterization of certain locally rank-one
symmetric spaces', {\it J. Differential Geom.} {\bf 28} (1988),
187--202.

\bibitem{Fie} Fiedler, B.,
`Determination of the structure of algebraic curvature tensors by means of
Young symmetrizers', Seminaire Lotharingien de Combinatoire {\bf B48d}
(2003a). 20 pp. Electronically published: http://www.mat.univie.ac.at/$\sim$slc/;
see also math.CO/0212278.

\bibitem{E-K-R} Garc{\'\i}a-R{\'\i}o, E., Kupeli, D. N., and V{\'a}zquez-Lorenzo, R., {\bf
Osserman Manifolds in Semi-Riemannian Geometry}, Lecture notes in
Mathematics, Springer-Verlag.

\bibitem{Gilkey-Swann-Vanhecke} Gilkey, P., Swann, A., and Vanhecke, L., `Isoparametric
geodesic spheres and a conjecture of Osserman concerning the
Jacobi operator', {\it Quart. J. Math. Oxford} {\bf 46} (1995), 299--320.

\bibitem{gilkey-p-Osserman} Gilkey, P., `Algebraic curvature
tensors which are $p$-Osserman', {\it Differential Geom. Appl.}
{\bf 14} (2001), 297-311.

\bibitem{G01} Gilkey, P., `Bundles over projective spaces and algebraic curvature
tensors', {\it J. Geom} {\bf 71} (2001), 54--67.

\bibitem{Gilkey-book} Gilkey, P., {\bf Geometric Properties of Natural Operators
Defined by the Riemann Curvature Tensor}, World Scientific Press,
Singapore (2001).

\bibitem{GHS} Glover, H., Homer,  W., and Stong, R., `Splitting the tangent bundle of
projective space', {\it Indiana Univ. Math. J.} {\bf 31} (1982), 161--166.

\bibitem{nikolayevsky-1} Nikolayevsky, Y., `Osserman manifolds of dimension $8$', {\it
Manuscr. Math.} {\bf 115} (2004), 31--53.

\bibitem{nikolayevsky-2} Nikolayevsky, Y.,  `Osserman Conjecture in dimension $n \ne 8,
16$', {\it Math. Ann.} {\bf 331} (2005), 505-522.

\bibitem{Osserman} Osserman, R.,
  `Curvature in the eighties',
  {\it Amer. Math. Monthly}, {\bf 97} (1990), 731--756.

\bibitem{StanilovVidev} Stanilov, G., and Videv,  V., `On a generalization of the
Jacobi operator in the Riemannian geometry',  {\it Annuaire Univ. Sofia Fac. Math.
Inform.} {\bf 86} (1992), 27--34.

\end{thebibliography}
\end{document}